\documentclass[10pt]{amsart} 
\usepackage{pslatex}
\theoremstyle{plain}
\newtheorem{theo}{Theorem}
\newtheorem{lemma}[theo]{Lemma}

\newtheorem{Definition}[theo]{Definition}
\newtheorem{Remark}[theo]{Remark}{}

\usepackage{pstricks}
\usepackage{amscd}
\usepackage{amssymb}
\usepackage{amsmath}

\begin{document}
\title{ Markov Towers and Stochastic Properties of Billiards}

\author{Domokos Sz\'asz,  Tam\'as Varj\'u}
\address{Budapest University of Technology,
Mathematical Intitute  and Center for Applied Mathematics,
Budapest, Egry J. u. 1, Hungary H-1111}
\email{szasz@math.bme.hu, kanya@math.bme.hu}
\thanks{Research supported
 by the Hungarian National Foundation for Scientific Research grants
  No.\ T32022 T26176 and Ts040719, and by FKFP 0058/2001}
\thanks{{\it 2000 Mathematics Subject Classification} {\rm Primary: 37D50, 60F05; Secondary: 37A50, 37A60}}
\thanks{{\it Key Words and Phrases:} {\rm Markov towers, local limit theorem, Sinai billiard, Lorentz process, recurrence}}
\dedicatory{Dedicated to Anatole Katok on the occasion of his 60th birthday.}
\maketitle

\begin{abstract} 
Markov partitions work most efficiently for Anosov systems or for Axiom A
systems. However, for hyperbolic dynamical systems which are either singular 
or whose hyperbolicity is nonuniform, the construction of a Markov
partition, which in these cases is necessarily countable, is a rather 
delicate issue even when such a construction exists. An additional problem
is the use of a countable Markov partition for proving probabilistic 
statements. For a wide class of hyperbolic systems, L. S. Young, in 1998, 
constructed so called Markov towers, which she could apply successfully
to establish nice, for instance, exponential correlation decay, and, 
moreover, as a consequence, a central limit theorem.
The aim of this survey is twofold. First we show how the Markov
tower construction is applicable for obtaining finer stochastic properties,
like a local limit theorem of probability theory. Here the fundamental method 
is the study of the spectrum of the Fourier transform of the Perron--Frobenius
operator. These ideas and results are applicable to all systems Young
has been considering. Second, we survey the problem of recurrence
of the planar Lorentz process. As an application of the results from
the first part, we obtain a dynamical proof of recurrence for the finite 
horizon case. Here basically different proofs were given by K. Schmidt, in
1998, and J.-P. Conze, in 1999. As another application we can also treat
the infinite horizon case, where already the global limit theorem is
absolutely novel. It is not a central one, the scaling is $\sqrt {n \log n}$
in contrast to the classical $\sqrt n$ one. Beyond thus giving a rigorous
proof for earlier heuristic ideas of P. Bleher, which used three 
delicate and hard hypotheses, we can  also a) verify the local version
of this limit theorem for the free flight function and b) prove the 
recurrence of the planar Lorentz process in the infinite horizon case.
\end{abstract}

\section{Introduction}
\label{sec:intro}

Since --- following some ideas of Hadamard --- M. Morse introduced the 
concept of symbolic dynamics, the method got more and more 
extensively used to study topological, and later also ergodic and
stochastic, properties of dynamical systems possessing some hyperbolic 
behaviour. 

On the one hand, ``the idea of coding and semiconjugacies with 
topological Markov chains yields remarkably precise results concerning
topological entropy, the growth of periodic orbits, the presence
of orbits of various periods, and the structure of maps with zero 
topological entropy'' \cite{KH}.
On the other hand, through the achievements of Bowen (cf.\ \cite{B 75},
Ruelle (cf.\ \cite{R 78} and of Sinai (most notably his work \cite{S 72}
relating symbolic dynamics and Gibbs states of statistical physics),
almost invertible semiconjugacies and conjugacies provided by 
Markov partitions made it possible to demonstrate exponential correlation decay
and further nice and useful stochastic properties for Axiom A systems ---
 and later for 
more general ones, too. Thus, for a long time it, quite naturally, seemed so
that the construction of Markov partitions is `the method' for obtaining
effective statistical statements for more complicated systems as well.

However, in course of the work of the Moscow school on billiards and of Benedicks and
Young on the H\'enon map it became clear that a Markov partition is a
too delicate construction if one wants to relax the assumptions of
smoothness or of the uniform hyperbolicity of the maps in question.
In fact, a warning might have come earlier from Bowen's result
showing that, even in the nicest systems,
in the multidimensional ($d >2$) 
case the boundaries of the elements of any Markov partition behave wildly,
in particular, they
are not smooth \cite{B 78}. In addition,
for instance,  in discontinuous systems like billiards, the
local invariant manifolds are, indeed, arbitrarily short, and as
a consequence the Markov partition is necessarily countable and its elements
are products of Cantor sets. Then to adapt the boundaries of a possible Markov
partition in a Markov way is an extremely delicate issue even in the case
when such a construction was successfully established (cf.\ \cite{BS 80}).

It is not our aim to go into more details here since there exists a
quite recent and excellent 
survey \cite{ChY}, which, on the one hand, gives a comprehensive historical
overview, and, on the other hand, explains the way out: the `weaker' 
construction of a Markov tower. This construction was designed  by L. S. Young
 \cite{Young}
and it works for a wide class of systems with some hyperbolicity, among others
for Anosov and Axiom A systems, two dimensional hyperbolic systems with 
singularities, two dimensional Sinai billiards with a finite horizon, 
hyperbolic unimodal maps and hyperbolic  H\'enon maps.

The construction of a Markov
partition and the resulting symbolic dynamics opened the way in a 
straightforward manner 
to put the probabilistic arsenal of --- appropriately mixing --- 
stationary stochastic processes into action.
In the case of a Markov tower this
connection is not straightforward, and our actual aim is, indeed, to
understand and to  discuss
the probabilistic approach in the case of a Markov tower.
We note that, as we will see in subsection 2.2, the tower also leads to
a (countable) Markov  partition, but
its properties, apart from its formal ones, are 
quite different from those of a traditional
Markov partition. In particlular, it does not seem to provide a 
flexible Markov approximation. This is why in its applications
new methods are needed and, in fact, their discussion is our main aim
here.

The paper is organized as follows: in section \ref{sec:mar} we recall 
the axioms of systems, which we are going to deal with, and 
briefly describe the tower construction. 
In section 3 we analyze how one can establish stochastic properties,
in general, and further  finer stochastic properties, like
local limit theorems, in particular. The choice of local limit theorems may seem 
eventual but it is not. This will be clear from section 4, where we
apply the local limit theorems to planar dispersing billiards.
In doing so, beside attaining local theorems for them, we also 
obtain
\begin{itemize}
\item a dynamical proof of recurrence for the planar Lorentz process
with a finite horizon (in fact, partly abstract ergodic-theoretic
proofs were given by K. Schmidt \cite{Sch} and Conze \cite{Conze}, which were, however, also
using the central limit theorem);
\item the first rigorous proof for a noncentral limit theorem for the 
displacements of the planar Lorentz process with an infinite horizon;
this result was conjectured by
 an earlier nonrigorous, heuristic argument of Bleher \cite{Bleher}
based on three hard hypotheses (which, in fact, still do not follow from our 
approach);
\item the first 
proof of recurrence for the planar Lorentz process with an infinite horizon.
\end{itemize}

\section{Markov tower}
\label{sec:mar}
\subsection{The Product Set}
\label{subsec:You}

The technique developed in \cite{Young} allows to handle stochastic
 properties of systems
\begin{enumerate}
  \item whose every power is ergodic;
  \item which satisfy several technical assumptions well-known from hyperbolic theory;
  \item whose phase space $X$ contains a subset $\Lambda$ with a hyperbolic product structure;
  \item where the return time into $\Lambda$ has an exponentially decaying tail.
\end{enumerate}
This class contains planar dispersing billiards with
 both bounded or  unbounded free flight (i.\ e.\ with 
a finite resp.\ an infinite horizon), logistic interval maps,
expanding maps with neutral fixed points, piecewise hyperbolic maps, H\'enon attractors, their generalisations, and
certain partially hyperbolic systems. The rest of this subsection is devoted to the precise definitions.

We start with describing precisely the models we are going to deal with.
Let $T$ be a $C^{1+\epsilon}$ diffeomorphism with singularities of a compact Riemannian manifold $X$ with boundary. More
precisely, there exists a finite or countably infinite number of pairwise disjoint open regions $\{X_i\}$ whose boundaries are $C^1$
submanifolds of codimension 1, and finite volume such that $\cup X_i = X$, $T\big|_{\cup X_i}$ is $1-1$ and $T\big|_{X_i}$
can be extended to a $C^{1+\epsilon}$-diffeomorphism of $\bar X_i$ onto its image. The
Riemannian measure will be denoted by $\mu$, and if $W\subset X$ is a submanifold, then $\mu_W$ will denote the induced
measure. The invariant Borel probability measure will be denoted by $\nu$.  \begin{Definition} An embedded disk
$\gamma\subset X$ is called an \emph{unstable manifold} or an \emph{unstable disk} if $\forall x,y\in\gamma, \enskip
d(T^{-n}x,T^{-n}y)\rightarrow 0$ exponentially fast as $n\rightarrow\infty$; it is called a \emph{stable manifold} or a
\emph{stable disk} if $\forall x,y\in\gamma, \enskip d(T^nx,T^ny)\rightarrow 0$ exponentially fast as
$n\rightarrow\infty$. We say that $\Gamma^u=\{\gamma^u\}$ is a \emph{continuous family of $C^1$ unstable disks} if the
following hold:
\begin{itemize}
  \item $K^s$ is an arbitrary compact set; $D^u$ is the unit disk of some $\mathbb{R}^n$;
  \item $\Phi^u\colon K^s\times D^u\rightarrow X$ is a map with the property that
  \begin{itemize}
    \item $\Phi^u$ maps $K^s\times D^u$ homeomorphically onto its image,
    \item $x\rightarrow\Phi^u\mid(\{x\}\times D^u)$ is a continuous map from $K^s$ into the space of $C^1$ embeddings of
      $D^u$ into X,
    \item $\gamma^u$, the image of each $\{x\}\times D^u$, is an unstable disk.
  \end{itemize}
\end{itemize}
\emph{Continuous families of $C^1$ stable disks} are defined similarly.
\end{Definition}

\begin{Definition}
  We say that $\Lambda\subset X$ has a \emph{hyperbolic product structure} if there exist a continuous family of
  unstable disks $\Gamma^u=\{\gamma^u\}$ and a continuous family of stable disks $\Gamma^s=\{\gamma^s\}$ such that
  \begin{enumerate}
    \item[(i)] $\dim \gamma^u+\dim\gamma^s=\dim X$
    \item[(ii)] the $\gamma^u$-disks are transversal to the $\gamma^s$-disks with the angles between them bounded away
      from $0$;
    \item[(iii)] each $\gamma^u$-disk meets each $\gamma^s$-disk in exactly one point;
    \item[(iv)] $\Lambda=(\cup\gamma^u)\cap(\cup\gamma^s)$.
  \end{enumerate}
\end{Definition}
\begin{Definition}
  Suppose $\Lambda$ has a hyperbolic product structure. Let $\Gamma^u$ and $\Gamma^s$ be the defining families for
  $\Lambda$. A subset $\Lambda_0\subset\Lambda$ is called an \emph{$s$-subset} if $\Lambda_0$ also has a hyperbolic
  product structure and its defining families can be chosen to be $\Gamma^u$ and $\Gamma^s_0$ with
  $\Gamma^s_0\subset\Gamma^s$; \emph{$u$-subsets} are defined analogously. For $x\in\Lambda$, let $\gamma^u(x)$ denote
  the element of $\Gamma^u$ containing $x$.
\end{Definition}

\begin{Definition}
  We call $(X,T,\nu)$ a \emph{Young system}, if the following Properties {\bf (P1)--(P8)} are true:
\end{Definition}

\begin{enumerate}
\item[\bf (P1)] There exists a $\Lambda\subset X$ with a hyperbolic product structure and with
  $\mu_\gamma\{ \gamma\cap\Lambda \}>0$ for every $\gamma\in\Gamma^u$.
\item[\bf (P2)] There is a countable number of disjoint $s$-subsets $\Lambda_1,\Lambda_2,\dots\subset\Lambda$ such that
  \begin{itemize}
  \item on each $\gamma^u$-disk $\mu_{\gamma^u}\{(\Lambda\setminus\cup\Lambda_i)\cap\gamma^u\}=0$;
  \item for each $i$, $\exists R_i\in\mathbb{Z}^+$ such that $T^{R_i}\Lambda_i$ is a $u$-subset of $\Lambda$;
  \item for each $n$ there are at most finitely many $i$'s with $R_i=n$;
  \item $\min R_i\geq$ some $ R_0$ depending only on $T$
  \end{itemize}
\item[\bf (P3)] For every pair $x,y\in\Lambda$, we have a notion of \emph{separation time} denoted by $s_0(x,y)$. If
  $s_0(x,y)=n$, then the orbits of $x$ and $y$ are thought of as being ``indistinguishable'' or ``together'' through their
  $n^\mathrm{th}$ iterates, while $T^{n+1}x$ and $T^{n+1}y$ are thought of as having been ``separated.''
  (This could mean that the points have moved a certain distance apart, or have landed on opposite sides of a
  discontinuity manifold, or that their derivatives have ceased to be comparable.) We assume:
 \begin{enumerate}
   \item[(i)] $s_0\geq 0$ and depends only on the $\gamma^s$-disks containing the two points;
   \item[(ii)] the number of ``distinguishable'' n-orbits starting from $\Lambda$ is finite for each $n$;
   \item[(iii)] for $x,y\in\Lambda_i, \enskip s_0(x,y)\geq R_i+s_0(T^{R_i}x,T^{R_i}y);$
 \end{enumerate}
\item[\bf (P4)] Contraction along $\gamma^s$ disks. There exist $C >0$ and $\alpha < 1$ such that for
  $y\in\gamma^s(x),\enskip d(T^nx, T^ny)\leq C\alpha^n\enskip\forall n\geq 0$.
\item[\bf (P5)] Backward contraction and distorsion along $\gamma^u$. For $y\in\gamma^u(x)$ and $0\leq k\leq
  n<s_0(x,y)$, we have
  \begin{enumerate}
  \item[(a)] $d(T^nx,T^ny)\leq C\alpha^{s_0(x,y)-n}$;
  \item[(b)] \[\log\prod_{i=k}^n\frac{\det DT^u(T^ix)}{\det DT^u(T^iy)}\leq C\alpha^{s_0(x,y)-n}.\]
  \end{enumerate}
\item[\bf (P6)] Convergence of $D(T^i|\gamma^u)$ and absolute continuity of $\Gamma^s$.
  \begin{enumerate}
  \item[(a)] for $y\in\gamma^s(x)$,\[\log\prod_{i=n}^\infty\frac{\det T^u(T^ix)}{\det T^u(T^iy)}\leq C\alpha^n\quad\forall
    n\geq 0.\]
  \item[(b)] for $\gamma,\gamma'\in\Gamma^u$, if $\Theta\colon\gamma\cap\Lambda\rightarrow\gamma'\cap\Lambda$ is defined
    by $\Theta(x)=\gamma^s(x)\cap\gamma'$, then $\Theta$ is absolutely continuous and
    \[\frac {d(\Theta_*^{-1}\mu_{\gamma'})} {d\mu_\gamma} (x) = \prod_{i=0}^\infty \frac {\det DT^u(T^ix)} {\det
    DT^u(T^i\Theta x)}.\] 
  \end{enumerate}
\item[\bf (P7)] $\exists C_0>0$ and $\theta_0<1$ such that for some
  $\gamma\in\Gamma^u$,\[\mu_\gamma\{x\in\gamma\cap\Lambda:R(x)>n\}\leq C_0\theta_0^n\quad \forall n\geq 0;\]
\item[\bf (P8)] $(T^n,\nu)$ is ergodic $\forall n\geq1$, where $\nu$ is
the SRB-measure corresponding to the system $(X, T)$.
\end{enumerate}

\begin{Remark}
Some explanation: Properties {\bf (P1-7)} do not involve the measure $\nu$.
Using them, in fact, Young constructs the SRB-measure $\nu$, so the system
satisfying {\bf (P1-7)} is a Young sysstem if for the so constructed 
SRB-measure  {\bf (P8)} holds, too.
\end{Remark}

\subsection{The Tower}

Now we will define the \emph{Markov extension}, the actual
\emph{ Markov tower}. 
Let $R:\ \Lambda \to \mathbb Z_+$ be the function which is
$R_i$ on $\Lambda_i$, and let
\[\Delta\stackrel{\mathrm{def}}{=}\{(x,l): x\in\Lambda;\enskip l=0,1,\dots,R(x)-1\}\] 
and define 
\[F(x,l)=\left\{
  \begin{array}{ll}
(x,l+1) & \text{if}\quad l+1<R(x)\\
(T^Rx,0) & \text{if}\quad l+1=R(x)
  \end{array}\right. \]
We will refer to $\Delta_l=\{(x,j)\mid (x,j) \in \Delta,\ j=l \}$ 
as the $l^\mathrm{th}$ level of the tower $\Delta$. Young also has a construction for $\tilde
\nu$, the SRB measure of the extension, for which the pushforward is $\nu$, and $J(F)\equiv 1$ except on
$F^{-1}(\Delta_0)$.
Thus, the Markov Tower is the dynamical system $(\Delta, F, \tilde \nu)$.

On the tower a {\emph Markov partition}
 $\mathcal{D}$ can be defined, with the following properties:
\begin{enumerate}
\item[(a)] $\mathcal{D}$ is a refinement of the partition $\Delta_l$.

\item[(b)]  By  denoting by $\mathcal{D}_l$ the partition
 $\mathcal{D}|\Delta_l$, 
$\mathcal{D}_l = \{\Delta_{l,j} \mid j=1, \dots, j_l\}$ 
has only a finite number of elements and
each one is the union of a collection of $\Lambda_i$'s;
\item[(c)] $\mathcal{D}_l$ is a refinement of $F\mathcal{D}_{l-1}$;
\item[(d)] if $x$ and $y$ belong to the same element of $\mathcal{D}_l$, then $s_0(F^{-l}x,F^{-l}y)\geq l$;
\item[(e)] if $R_i=R_j$ for some $i\not=j$, then $\Lambda_i$ and $\Lambda_j$ belong to different elements
  of $\mathcal{D}_{R_i-1}$.
\end{enumerate}
Let $\Delta^*_{l,j}=\Delta_{l,j}\cap F^{-1}(\Delta_0)$. We think of $\Delta_{l,j}\setminus\Delta^*_{l,j}$ as ``moving
upward'' under $F$, while $\Delta_{l,j}^*$ returns to the base.

\medbreak 
It is natural to \emph{redefine the separation time} to be $s(x,y)\stackrel{\mathrm{def}}{=}$ the
largest $n$ such that for all $i\leq n,\enskip F^ix$ and $F^iy$ lie in the same element of $\{\Delta_{l,j}\}$.  We claim
that \textbf{(P5)} is valid for $x,y\in\gamma^u\cap\Delta_{l,j}$ with $s$ in the place of $s_0$. To verify this, first
consider $x,y\in\Lambda$. We claim that $s(x,y)\leq s_0(x,y)$. If $x,y$ do not belong to the same $\Lambda_i$, then this
follows from rule (d) in the construction of $\mathcal{D}_l$; if $x,y\in\Lambda_i$, but $T^Rx,T^Ry$ are not contained in
the same $\Lambda_j$, then $s(x,y)=R_i+s(T^Rx,T^Ry)$, which is $\leq s_0(x,y)$ by property \textbf{(P3)},(iii) of $s_0$,
and so on. In general, for $x,y\in\Delta_{l,j}$, let $x_0=F^{-l}x,\enskip y_0=F^{-l}y$ be the unique inverse images of $x$
and $y$ in $\Delta_0$. Then by definition $s(x,y)=s(x_0,y_0)-l$, and what is said earlier on about $x_0$ and $y_0$ is
equally valid for $x$ and $y$.

\textit{From here on $s_0$ is replaced by $s$ and \textbf{(P5)} is modified accordingly.}

\section{Stochastical Properties}

\label{sec:sto}

As explained above, the tower construction actually provides a countable
Markov partition. Below we first remind the reader how Young exploits this  partition to obtain
stochastic behaviour. The first goal is to establish exponential decay of correlations, where, of course,
 Property {\bf (P7)} plays a decisive role. 

\subsection{The Perron--Frobenius Operator and the Doeblin--Fortet Property}

\paragraph{1.} The starting point is to \emph{factorise the dynamics} by  a factorisation along stable 
manifolds of $\Delta$.\@ The advantage is that
this dynamics will behave as an expanding map, an appropriate
 object to study via the Perron--Frobenius operator. 
Let $\bar{\Delta}:=\Delta/\sim$ where $x\sim y$
iff $y\in\gamma^s(x).$ Since $F$ takes $\gamma^s$-leaves to $\gamma^s$-leaves, the quotient dynamical system
$\bar{F}\colon\bar{\Delta}\rightarrow\bar{\Delta}$ is clearly well defined.

The construction of
$\bar m$ is not trivial but is quite standard. Young obtains it following \cite{B 75}.
A simple property of  $\bar{m}$ is: let $\bar{m}|\bar{\Delta}_l$ be the measure induced from the natural
identification of $\bar{\Delta}_l$ with a subset of $\bar{\Delta}_0$, so that $J(\bar{F})\equiv 1$ except on
$\bar{F}^{-1}(\bar{\Delta}_0),$ where $J(\bar{F})=J(\overline{T^R}\circ\bar{F}^{-(R-1)}).$ 
We note that for the factorised map Young also proves a distorsion
property  with a weaker constant $\beta:\ 1 > \beta > \sqrt \alpha$. 

To investigate Birkhoff sums we have to associate a
function $\bar{f}:\bar{\Delta}\rightarrow\mathbb{R}$ to each observable
$f:M\rightarrow\mathbb{R}$. We can pull back $f$ to the tower, and find an
other function cohomologous to this one, which is constant along stable
manifolds. This method is described for example in \cite{PP}. 

\paragraph{2.} The analytic tool of investigation is the \emph{Perron--Frobenius operator}:
\[P(\bar{\varphi}(\bar{x}))=\sum_{\bar{y}:\bar{F}\bar{y}=\bar{x}} \frac {\bar{\varphi}(\bar{y})} {J\bar{F}(\bar{y})} .\]

\paragraph{3.} The technique is based upon the \emph{spectral properties of the Perron--Frobenius
operator}. For this purpose we have to introduce \emph{suitable Banach spaces}, where $P$ has
a nice spectrum. Actually for the method introduced by Doeblin and Fortet we
need a pair of function spaces ${\mathcal C}$ and ${\mathcal L}$ (usually with some
supremum-like and Lipschitz-like norms) such that ${\mathcal L}\leq{\mathcal C}$,
$\|\cdot\|_{\mathcal C}\leq\|\cdot\|_{\mathcal L}$, and the inclusion of ${\mathcal L}$
into ${\mathcal C}$ be a compact operator.

If we have such a pair of Banach spaces, and we can prove, that $\exists N,K$,
and $\tau<1$ \[ \|P^N\bar{\varphi}\|_{\mathcal L} \leq \tau
\|\bar{\varphi}\|_{\mathcal L} + K\|\bar{\varphi}\|_{\mathcal C}\] then by knowing
that $\forall i \enskip T^i$ is ergodic we have that the spectrum of $P$ on ${\mathcal L}$
is contained in a disk with radius strictly smaller than one, except that 1 is
an eigenvalue with multiplicity one, and the corresponding eigenfunction is
the invariant density (cf. \cite{IT-M}).

This kind of estimate captures the uniform expanding feature of the dynamics, or contraction of the $P$ operator. In order
to derive this so called \emph{Doeblin--Fortet property of the transfer operator} (in
the theory of dynamical systems often called the \emph{Lasota--Yorke
  property}, cf. \cite{LY}), Young uses an exponential factor in the function norms. $\bar{F}$ on the tower has Jacobian 1, when moving
upwards, and the tower is usually unboundedly high, so we have to pretend expanding at least in the norms:
\begin{align*}
  \|\bar{\varphi}\|_{\mathcal C} &:= \sup_{l,j}
  \left|\bar{\varphi}\big|_{\Delta_{l,j}}\right|_{\infty}
  e^{-l\epsilon}, \\
  \|\bar{\varphi}\|_{{\mathcal L}} &:= \|\bar{\varphi}\|_{\mathcal C} + \left(
    {}\underset {\scriptsize \bar{x},\bar{y}\in\bar{\Delta}_{l,j}}
    {\mathrm{ess \enskip sup \enskip}} \frac
    {|\bar{\varphi}(\bar{x})-\bar{\varphi}(\bar{y})|}
    {\beta^{s(\bar{x},\bar{y})}} \right) e^{-\epsilon l}.
\end{align*}
The denominator in the second definition is a natural distance for the points on the tower, so it is really a
Lipschitz-like norm. The aforementioned trick is in the $e^{-l\epsilon}$ term. This allows that a function on
$\bar{\Delta}$ which is exponentially increasing with the height of the tower, to be in ${\mathcal C}$, if this growth is
moderate. This also means, that $P$ is contracting the norm in the middle of the tower, where the Jacobian is 1. The
constant $\epsilon$ should be chosen carefully in order to hold back enough contraction when  a Markov
return occurs. Roughly speaking $\epsilon$ has to be smaller than the smallest positive Lyapunov exponent.

\subsection{The Central Limit Theorem}

The aforementioned spectral picture is essentially equivalent to the exponential decay of
correlations for function pairs in ${\mathcal L}$. As a matter of fact, the boundedness 
of the functions in question is also needed. 
After the argument presented in \cite{Young} this leads to the exponential correlation
decay of bounded, piecewise H\"older functions on $M$.
For the same class of functions Young immediately 
gets the \emph{Central Limit Theorem} (CLT) by
checking the conditions of a theorem by Keller \cite{Ke}.  
Here we make a simple but important clarification.
The traditional and by far the most widely used method for establishing the CLT
is to apply Fourier transforms. Keller \cite{Ke} (and thus also \cite{Young}) can elude this
by referring to a nice and useful theorem of Gordin \cite{G}. (As a matter 
of fact,
Fourier transform are, indirectly, still applied, since Gordin 
constructed a martingale approximation and
used the martingale CLT. But in proving the CLT for martingales again 
Fourier transform is the
method! Moreover, the error term in the martingale approximation
is so large that it obviously excludes the applicability
of this approach in proofs of finer statements, for instance, in 
those of local limit theorems.)

Let we recall two basic results from \cite{Young}:
the first one on the exponential decay of correlations and 
the second one on the CLT. Notation:
\[
\mathcal H_\eta = \left \{\varphi: M \to \mathbb R \bigm| \exists A >0
\ \text{such that for}\
\forall \ x,y\in M\ |\varphi(x) - \varphi(y)| \le Ad(x,y)^\eta \right\}.
\]

\begin{theo}
(\cite{Young}) For any $\eta$, there exists $\tau < 1$ such that for all
$\varphi, \psi \in \mathcal H_\eta$ there exists a
$C=C(\varphi, \psi)$ such that
\[
\left | \int (\varphi \circ T^n) \psi d\nu - \int \varphi d\nu
\int \psi d\nu \right | \le C \tau^n \qquad \forall n \ge 1
\]
\end{theo}

\begin{theo}(\cite{Young})
Assume $\varphi\in \mathcal H_\eta$ and $\int \varphi d\nu=0$. Then
\[
\frac{1}{\sqrt n} \sum_{i=0}^{n-1} \varphi \circ T^i 
\Longrightarrow \mathcal N(0, \sigma^2)
\]
for some $\sigma \ge 0$. Moreover, $\sigma = 0$ iff $ \varphi = \psi \circ T
- \psi$ for some $\psi \in L^2(\nu)$. (Here $\Longrightarrow$
denotes weak convergence of probability distributions and $\mathcal 
N(m, \sigma^2)$ 
denotes the normal distribution with mean $m$ and variance $\sigma^2$.)
\end{theo}

For simplicity, we have formulated these results for one-dimensional
random variables, and their extension to vector valued functions
is straightforward.

\subsection{Local limit theorem}
\label{subsec:llt}

For illustrating a local CLT as compared to the widely used (global) CLT, take a simple symmetric 
random walk (SSRW) on $\mathbb Z^d$.
 So let
 $W_n = X_1 + \dots X_n$, where $X_1, \dots, X_n, \dots$
are independent, identically distributed random variables with the common distribution
$P(X_i = \pm e_j) = \frac{1}{2d};\ 1 \le j \le d$ for all $i \in 
\mathbb Z_+$ (here the $e_j$s are the standard unit vectors of $\mathbb Z^d$).
 Then, of course, 
the CLT says that $P({W_n} \in {\sqrt n} A) \to \Phi (A)$ as $n \to \infty$, where $\Phi(A) = \int_A
\phi(s)ds$ and $\phi(s) = (2\pi)^{-d/2} \exp{(-\frac{s^2}{2})}$ is the 
$d$-dimensional Gaussian density.
In other words it describes the asymptotics of a sequence of sets increasing like $\sqrt n$.
In contrast, for the SSRW  the local CLT (LCLT) 
says that $n^{d/2} P(W_n = [s\sqrt n]) \to \phi(s)$ as $n \to 
\infty$, i.\ e.\ it describes the asymptotics of a sequence of sets of fixed size,
consequently it is, indeed, local! 

For stating our main theorem we have to fix some notations first. For a fixed $f:X \to \mathbb{R}^d$ denote the average
$\int f d\nu = a$, and \[ S_n^f(x) 
= \sum _{k=0}^ {n-1} f (T^{k}x) \] the Birkhoff sum. Consider the smallest translated closed
subgroup $V+r\subseteq\mathbb{R}^d$ which supports the values of $f$ ($V$ is the group and $r$ is the translation). By
ergodicity of all powers of $T$, the support of $S_n$ is $V+nr$.
\begin{theo}\label{thm:main}(\cite{Main}
  Suppose that
  \begin{enumerate}
    \item $(X, T, \nu)$ is a Young system (cf.\ subsection \ref{subsec:You});
    \item $f$ is minimal: i.\ e.\ it is not cohomologous to a function for which the support in the above sense is
      strictly smaller.
    \item $f$ is nondegenerate: i.\ e.\  $\text{span} \left<V\right>=\mathbb{R}^d$, and
    \item $f$ is bounded and H\"older continuous.
  \end{enumerate}  
  Let $k_n\in V+nr$ be such that $\frac{k_n-na}{\sqrt{n}}\rightarrow k$.  Denote the distribution of $S^f_n-k_n$ by
  $\upsilon_n$. Then \[ n^{d/2} \upsilon_n\rightarrow \phi(k) l\] where $\phi$ is a nondegenerate normal density
  function with zero expectation, and $l$ is the uniform measure on $V$: product of counting measures and Lebesgue
  measures. The convergence is meant in the weak topology.
\end{theo} 
\begin{Remark}
  For nonminimal functions we can obtain an analogous result. The limit measure on the right hand side in this case is
  not necessarily uniform.
\end{Remark}

We know from classical analysis or from probability theory that the local behaviour of 
distributions (densities, measures,\dots) is connected to the tail behaviour of
the corresponding Fourier transforms. Therefore,
in the next subsection, we are going to study the Fourier 
transform of the Perron--Frobenius operator. 

\subsection{The Fourier Transform of the Perron--Frobenius Operator}

When one wants to obtain finer results than the CLT, then
Fourier transforms seem 
inevitable. Thus, in our setup, we have to define the
Fourier transform of the Perron--Frobenius operator: \[
P_t(\bar{\varphi}):=P(e^{it\bar{f}}\bar{\varphi})\] where $f$ is the
function for which the limit theorem is stated. Note that $P_0=P$!
Also, it is worth noting that $P_t^n \mathbf 1 = P^n(\exp{(it S_n^{\bar f})})$
and, moreover, $\mathbb EP_t^n \mathbf 1=\mathbb EP^n(\exp{(it S_n^{\bar f})})
 = \mathbb \exp{(itS_n^{\bar f})}$, the usual
characteristic function of $S_n^{\bar f}$ 
where $S_n^{\bar f}$ denotes the Birkhoff sum for the function $\bar f$.

The heart of this method is to expand the leading eigenvalue of $P_t$
--- analogously to the Taylor expansions around $0$  of characterisctic
functions of probability theory. Before
that, however, 
 we have to ensure its existence. We proved in \cite{Main} that $t\mapsto P_t$ is
continuous in the ${\mathcal L}$-norm, and by the stability
of the spectrum, for small values
of $t$ there exists $\lambda_t$, the perturbed value of 1 as an eigenvalue
with multiplicity one.  By proving the Doeblin--Fortet inequality for $P_t$ we get
that the rest of the spectrum will lie in a disk, with radius smaller than 1,
so it will not bother $\lambda_t$ to be the leading eigenvalue.

If the function $f$ is bounded and piecewise H\"older continuous on $M$, then we get the
second order Taylor expansion for the Fourier transform: \[P_t(\bar{\varphi})
= P( \bar{\varphi}) + it P\left( \bar{f} \bar{\varphi} \right) - \frac{t^2}2 P
\left( \bar{f}^2 \bar{\varphi}\right) + o(t^2) \left\| \bar{f}^2 \bar{\varphi}
\right\|_\mathcal{L}. \] From the assumptions it follows that $ \left\|
  \bar{f}^2 \bar{\varphi} \right\|_\mathcal{L}$ is finite. By an argument
presented in various forms in
\cite{Nag}, \cite{KSz2} and \cite{GH} this leads to the expansion of $\lambda_t$.

The philosophy explained above can already be combined with the classical 
proof of the local CLT sketched in the Appendix. Indeed, term $II$ is the same, term $I$
corresponds to the CLT. Term $III$ can be handled by using the ideas 
outlined above, while for bounding term $IIII$ we have used some compactness
arguments borrowed from \cite{AD}.

\section{The Planar Lorentz Process}
\label{sec:rec}

\subsection{Semidispersing billiards and Lorentz process}

In this subsection we summarize some basic properties of semidispersing billiards. Our aim is to introduce the most
important concepts and fix the notation which is essentially borrowed from \cite{KSSz}. Semidispersing billiards are
more or less hyperbolic dynamical systems with singularities.  Pesin's theory was extended to  these systems
in \cite{KS}.

A billiard is a dynamical system describing the motion of a point particle in a connected, compact domain $Q \subset
\mathbb{T}^d$. The boundary of the domain in assumed to be piecewise $C^3$-smooth. Inside $Q$ the motion is uniform while
the reflection at the boundary $\partial Q$ is elastic. As the absolute value of the velocity is a first integral of
motion, the phase space of the billiard flow is fixed as $M=Q\times S^{d-1}$ -- in other words, every phase point $x$ is
of the form $x=(q,v)$ with $q\in Q$ and $v\in \mathbb{R}^d,\ |v|=1$.  The Liouville probability measure $\mu$ on $M$ is
essentially the product of the Lebesgue measures, i.\ e.\ $d\mu= {\rm const.}\, dq dv$. The resulting dynamical system
$(M, S^{\mathbb{R}} , \mu)$ is the (toric) \emph{billiard flow}.

Let $n(q)$ denote the unit normal vector of a smooth component of the boundary $\partial Q$ at the point $q$, directed
inwards $Q$.  Throughout the paper we restrict our attention to \emph{semidispersing billiards}: we require for every
$q\in \partial Q$ the second fundamental form $K(q)$ of the boundary component to be nonnegative.

The boundary $\partial Q$ defines a natural cross section for the billiard flow. Namely consider \[ \partial M = \{ (q,v)
\mid q\in \partial Q,\left< v,n(q)\right> \ge 0 \}.\] This set actually has a natural bundle structure (cf.\
\cite{4geom}). The Poincar\'e section map $T$, also called the \emph{billiard map} is defined as the first return map
on $\partial M$. The invariant measure for the map is denoted by $\mu_1$, and we have $d\mu_1= {\rm const.} \left|\left<
    v,n(q)\right>\right|dqdv$. Throughout the paper we work with this discrete time dynamical system $(\partial
M,T,\mu_1)$. Recall the usual notation: for $(q, v) \in M$ one denotes
$\pi(q, v) = q$ the natural projection.

The \emph{Lorentz process} is the natural $\mathbb{Z}^d$ cover of a toric billiard. More precisely: consider
$\Pi:\mathbb{R}^d \to \mathbb{T}^d$ the factorisation by $\mathbb{Z}^d$. Its fundamental domain $D$ is a $d$-dimensional
cube (semiopen, semiclosed) in $\mathbb{R}^d$, so $\mathbb{R}^d = \cup_{z \in \mathbb{Z}^d} (D+z)$, where $D+z$ is the
translated fundamental domain.

By denoting $\tilde Q = \Pi^{-1} Q$, $\tilde M = \tilde Q \times S^{d-1}$, etc., the Lorentz dynamics is $(\tilde M,
\{\tilde {S}^t\mid t \in \mathbb{R}\}, \tilde \mu)$ and its Poincar\'e section map is ($\partial \tilde M, \tilde T,
\tilde \mu_1)$.  The \emph{free flight function} $\tilde\psi: \partial \tilde M \to \mathbb{R}^d$ is defined as follows:
$\tilde \psi(\tilde x)=\tilde q(T\tilde x)-\tilde q (\tilde x)$. The \emph{discrete free flight function} $\tilde \kappa :
\partial \tilde M \to \mathbb{Z}^d$ is defined as follows: $ \tilde \kappa (\tilde x) = \iota (\tilde T \tilde x) - \iota
(\tilde x)$, where $ \iota (\tilde x) = z $ if $ \tilde x \in Dz$. Observe finally, that $\tilde\psi$ and $\tilde \kappa$
are invariant under the $\mathbb{Z}^d$ action, so there are $\psi$ and $\kappa$ functions defined on $\partial M$, such
that $\tilde\psi= \Pi^*\psi$ and $\tilde \kappa = \Pi^* \kappa$. Actually for our purposes it will be more convenient to
choose the fundamental domain in such a way that $\partial \tilde Q \cap \partial D=\emptyset$. In this way $\kappa$ will
be continuous.

\subsection{LCLT and recurrence for the $d=2$ finite horizon case}

Consider the Lorentz process starting from the fundamental cell, a fixed isomorphic version of the 
fundamental domain. In the domain the starting phase point is random, it is
distributed according to the invariant measure of the corresponding torus-billiard. The relative position of the Lorentz
particle after the $n^{\mathrm{th}}$ collision is $S_n^\psi$, the discrete position is $S_n^\kappa$. So the event $A_n =
\left( S_n^\kappa = 0 \right)$ means that after the $n^{\mathrm{th}}$ collision the particle is again in the fundamental
domain. Recurrence means that it happens almost surely. 

\begin{theo} \cite{Sch}, \cite{Conze}, \cite{Main}
 The planar Lorentz process with a finite horizon is recurrent.
\end{theo}

For the proof we will use a stronger version of the well-known
Borel--Cantelli lemma. This version is due to Lamperti: 
\begin{lemma}(\cite{S})
  If for the sequence of events $\{A_n\}$
\begin{gather*}
  \sum_{k=1}^\infty \nu (A_k) = \infty \intertext{and some asymptotic independence holds:}
\liminf_{n\rightarrow\infty}
  \frac {\sum\limits _{j,k=1}^n \nu(A_jA_k)} {\left( \sum\limits_{k=1}^n \nu(A_k) \right)^2} <c \intertext{then $A_n$
    happens infinite often almost surely:} \nu \left( (q,v)\in \partial M \Biggm| \exists n_k \rightarrow \infty \enskip
    (q,v)\in \bigcap_{k=0}^\infty A_{n_k} \right) = 1
\end{gather*}
\end{lemma}

To check the first condition we have to deal with the asymptotic probabilities of the events $S_n^\kappa=0$. This is
exactly the region of the local limit theorem. Restrict ourselves to the $d=2$ finite horizon case! This case is known to
be a Young system i.\ e.\ it satisfies \textbf{(P1)--(P8)}, so the first condition of our LCLT is satisfied. For the second
condition we need to check the minimality of $\kappa$:

\begin{theo}(\cite{Main})
  $\kappa$ is minimal in the class of $\psi$.
\end{theo}

The proof is quite involved. Surprisingly it is related to arguments in \cite{BSCH91},
\cite{BSp 96} and \cite{BinSz} for establishing
the nondegeneracy of the covariance matrix in the CLT. To prove minimality, for each sublattice
$\mathbb{L}\subset\mathbb{Z}^2$ of finite index $n=\mathbb{L} : \mathbb{Z}^2$, we needed a periodic point
$n|\mathrm{per}(x)$ such that the Birkhoff sum $S^\kappa_{\mathrm{per}(x)} (x) \not \in \mathbb{L}$. To find this point we
used again Markov properties of $\Lambda$. Details can be found in \cite{Main}.

Since $\kappa$ is minimal and the values span the plane it is nondegenerate. Since it is piecewise constant it is also
H\"older continuous, so the LCLT applies. This gives $\nu ( S^\kappa_n=0 ) \sim \frac {\mathrm{const}} {n}$, and the sum
clearly diverges.

The second condition of Lamperti's Borel--Cantelli lemma contains intersection of recurrence events, so we also had to
prove a LCLT for joint distributions.

\section{Infinite Horizon}

Chernov \cite{Ch:inf} observed that the planar infinite horizon Sinai billiard also satisfies \textbf{(P1)--(P8)}, so
according to theorem \ref{thm:main} the LCLT also holds exactly the same manner as stated for the finite horizon case. 
Nota
bene: for bounded functions $f:M\rightarrow\mathbb{R}^d$ satisfying the conditions in Theorem 
\ref{thm:main} the asymptotics of the
Birkhoff sum $\nu ( S_n^f=0 ) \sim \frac c {n^{d/2}}$.

But the $\psi$ or the $\kappa$ free flight functions are not covered by this theorem since here $\kappa$ is not bounded,
$\psi$ is even not H\"older. This is not a surprise, since in this case a new phenomenon shows up. Former heuristic
arguments by Bleher \cite{Bleher} already suggested to expect that the moving particle is superdiffusive namely $\frac
{S_n^\psi} {\sqrt{n\log n}}$ will have a limit distribution. The reason is the following:

The infinite horizon condition is equivalent to the existence of collision free orbits in the phase space $M$. These
orbits form corridors in $\mathbb{R}^2$. Large free flights can occur by ``crossing'' one of these corridors. The smaller
the angle with the direction of the corridor is, the longer is the free flight.

\noindent\setlength{\unitlength}{\textwidth}
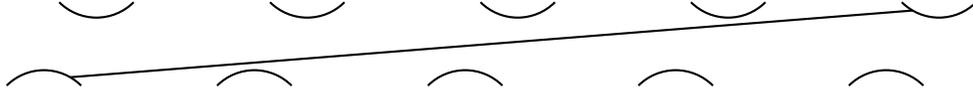
\begin{picture}(1,.15)
\psset{unit=\textwidth}

\multirput(.1,.15)(.2,0){5}{\psarc(0,0){.05}{225}{315}}
\multirput(.05,0)(.2,0){5}{\psarc(0,0){.05}{45}{135}}
\psline (.075,.0433)(.875,.1067)

\end{picture}

Computing the invariant measure for the small angle sets gives the asymptotics $\nu (|\kappa|=u) \sim \frac c {u^3}$. This
means that $\int |\kappa|^2 d\nu = \infty$, but any power with smaller exponent is integrable. As a matter of fact, the
distribution of $\kappa$ is in the nonnormal domain of attraction of the normal law, and its Fourier transform is: \[
\hat \kappa \stackrel {\mathrm{def}} {=} \int e^{it\kappa} d\nu = 1 + c |t|^2 \log |t| + O(|t|^2) \] (where $c$ is a
constant matrix) which means that if we would add independent copies of the same distribution then $\frac {S_n^*}
{\sqrt{n\log n}}$ would tend to a gaussian law. The $\log$ factor comes from the fact that the truncated $\kappa^x =
\kappa \mathbf{1}_{|\kappa| \leq x}$ has a variance of order $\log x$.

Direct geometrical calculations show that, when $\kappa$ is large, then the order of $\kappa\circ T$ is between
$\sqrt{\kappa}$ and $\kappa^2$. (There is a case when the trajectory hits the neighbouring scatterer on the same side of
the corridor before ``crossing''.\newline
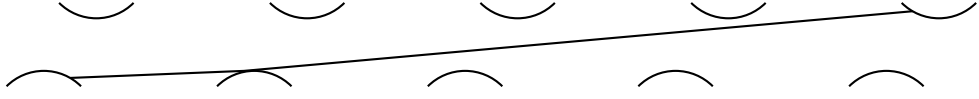
\begin{picture}(1,.15)
\psset{unit=\textwidth}

\multirput(.1,.15)(.2,0){5}{\psarc(0,0){.05}{225}{315}}
\multirput(.05,0)(.2,0){5}{\psarc(0,0){.05}{45}{135}}
\psline (.075,.0433)(.24,.05)(.875,.1067)

\end{picture}
In this case we change the Poincar\'e section and consider the sum of the small and the necessarily large free flight
vector as $\kappa\circ T$.) An important observation is that 
typically the next free flight will be in the regime of $\sqrt{\kappa}$. More precisely, for
any $\delta > \frac 12$ the probability $\nu ( \kappa\circ T > u^\delta \mid \kappa=u )\rightarrow 0$ as $u \rightarrow
\infty$. Even more the conditional expectation $\mathbb{E}_\nu (\kappa\circ T \mid \kappa)$ is of order $\sqrt \kappa$.
This means that though $|\kappa|^2$ is not integrable,
nevertheless, $\int \kappa (\kappa\circ T) d\nu <\infty$. The hope that the
autocorrelation may have a fast decay has lead to the conjecture that $S_n^\kappa$ asymptotically behaves the same way as
the sum of independent copies of $\kappa$.

\subsection{Limit theorems: global and local,  and recurrence for the $d=2$ infinite horizon case}

To reach the aforementioned limit theorem, and moreover the local limit theorem for $\kappa$ we used the symbolic space
$\bar \Delta$ as constructed in \cite{Young}, and investigated carefully how and where large values of $\bar \kappa$
appear on the tower. Clearly $\kappa$ is bounded on $\Lambda$. Since in one step it can grow at most to it is square, or
shrink to it is square root the time needed to reach the $\kappa=u$ set from $\Lambda$ is about $c\log\log u$, and before
returning to $\Lambda$ also the same amount of time is needed. By \textbf{(P7)} it is immediate that both the measure of
phase points, which spend $k$ iterates in the corridor before returning to the base of the tower, and both the measure of
phase points visiting the corridor $k$ times before returning to the base is exponentially small in $k$.

We had to replace the function norms in the definition of ${\mathcal C}$ and ${\mathcal L}$ to refer to $\bar \kappa$ on the
tower.  Instead of Young's $e^{-\epsilon l}$ factor we used \[\prod_{k=1}^l \min\left( e^{-\epsilon}, \bar
  \kappa^{-\delta} \circ \bar F^{-k} \right).\] Observe, that if $\kappa$ is bounded and $\delta$ is small enough this
gives back $e^{-\epsilon l}$. In that way we managed to achieve that $t \mapsto P_t$ be a countinuous mapping in both
function norms. This was not the case with the original norms, if the horizon is infinite. So $P_t$ can be considered as
a perturbation of $P$, and this also extends to the leading eigenvalue $\lambda_t$.

To obtain that the asymptotic behaviour is the same as in the case of the sum of independent copies, we needed to show
that \[ \lambda_t = 1 + c |t|^2 \log |t| + O(|t|^2) \] has the same kind of behaviour as the dynamically untouched $\hat
\kappa$. This is the key of the proof. While sums of independent copies give the product Fourier transform $\hat
\kappa^n$, that of the Birkhoff sums give $\int P_t^n (\mathbf{1}) d\nu$. The operator can be approximated by $\lambda_t^n
\bar \nu$ (here $\bar \nu$ is the leading
projection operator) up to an exponentially small error term coming from the rest of the
spectrum. Summarising: if $\hat \kappa (t)$ and $\lambda_t$ behave the same way as $t\rightarrow 0$, then the sum of
independent copies and the Birkhoff sum $S_n^\kappa$ behave the same way as $n\rightarrow\infty$.

There is a quite involved proof of the $\lambda_t$ expansion which is based on correlation estimates of powers of the
truncated $\bar \kappa$ and the eigenfunction related to $\lambda_t$. Details will appear in a technical paper.

\begin{theo}
  Suppose that the direction vectors of infinite collision free flights span the plane. Let $A \subset \mathbb{R}^2$. Then \[ \nu (S_n^\kappa \in \sqrt {n\log n} A) \rightarrow \int_A \phi(k)dk \] where
  $\phi$ is a nondegenerate gaussian density with zero expectation.
\end{theo}

\begin{Remark} The problem of the limiting behaviour of displacements in the case of an infinite horizon
has raised the interest of several people using different methods (very interesting works are
 \cite{Bleher} and \cite{ZE97}). 
It is worth mentioning that the computational
method of  \cite{ZE97} forecasts a non-Gaussian limit under the same scaling. .
\end{Remark}

\begin{theo}
  Suppose that the direction vectors of infinite collision free flights span the plane. Let $k_n\in\mathbb{Z}^2$ such that
  $\frac {k_n} {\sqrt{n\log n}}\rightarrow k$. Then \[ n \log n \enskip \nu (S_n^\kappa = k_n) \rightarrow \phi(k)\] where
  $\phi$ is a nondegenerate gaussian density with zero expectation.
\end{theo}

$\phi$ depends only on the corridor geometry. When computing the covariance one considers only the directions and widths
of corridors, and the bounding points of the corridors. This is a finite set of points on the scatterers, the geometry of
this finite set, and the curvature of the scatterers at these points are involved, but nothing else.  The sum of $\frac 1
{n\log n}$ diverges so the recurrence follows using an analogous argument as for the finite horizon case.  In the case
when all corridors are parallel one has to apply a nonisotropic scaling.

\begin{theo}
  Suppose, that all collision-free flights in the plane are parallel to the unit vector $w$. Consider the linear
  transformation $B_n$ which has the matrix $\scriptsize\left(
    \begin{array}{cc}
      \sqrt{n\log n} & 0 \\
      0 & \sqrt n
    \end{array}
  \right)$ in the basis $w,w^\perp$. Let $k_n\in\mathbb{Z}^2$ such that $B^{-1}_nk_n\rightarrow k$. Then \[ \det B_n
  \enskip \nu (S_n^\kappa = k_n) = n \sqrt {\log n} \enskip \nu (S_n^\kappa = k_n) \rightarrow \phi(k)\] where $\phi$ is a
  nondegenerate gaussian density with zero expectation.
\end{theo}

The sum of $\frac 1 {n\sqrt{\log n}}$ is also divergent, so recurrence is also obtained in this case.

\begin{Remark}
It is worth noting that in a recent manuscript of Gou\"ezel \cite{Gou} a related problem was investigated for
$1-D$ piecewise expanding maps with a neutral fixed point. The behaviour of billiard orbits in corridors is analogous to that 
of orbits of these $1-D$ maps near the neutral fixed point. Essential difficulties in our case arise from a) the larger 
dimension of the space; b) the not quite explicit form of the billiard dynamics; and c) emphatically from the fact that
in our case the function of interest is unbounded with a quite long tail whereas in \cite{Gou} it is bounded. His setup 
is more general since he is also considering stable limit laws in general, like \cite{AD}. The restriction of our
interest to the nonnormal domain of attraction of the Gaussian law came from the fact that our main concern was the free flight
function.
\end{Remark}

{\bf Acknowledgement.} The authors thank the refereee for his valuable remarks and suggestions.

\pagebreak
\section*{Appendix: A Classical Local CLT}

Here we recall, in the simplest case, the classical proof of Gnedenko \cite{Gn} of a local CLT.

For simplicity,
consider the case $d=1$. Following the notations of subsection \ref{subsec:llt}, our goal here is to
prove that, as $n \to \infty$,
$$
\sqrt n P(W_n=k_n) \to \phi(k)
$$
if $\frac{k_n}{\sqrt n} \to k$. Heuristically one expects that
$$
\sqrt n P(W_n=k_n) = \sqrt n P(\frac{W_n}{\sqrt n} = \frac{k_n}{\sqrt n}) 
\approx \sqrt n \phi(\frac{k_n}{\sqrt n}) \frac{2}{\sqrt n}
= \frac{1}{\pi}\int \exp(-is \frac{k_n}{\sqrt n}) \gamma(s) ds
$$
where $\gamma(s)=\exp(-\frac{s^2}{2})$ is the standard gaussian
characteristic function. Here we used the gaussian approximation 
suggested by the CLT and the 
Fourier inversion formula.

For a proof, let us turn to characteristic functions of $W_n$.
Denote by  $\xi(t)$ the common characteristic function of
the variables $X_i$. 
Then
$$
\sqrt n P(W_n=k_n)  = \sqrt n \frac{1}{2\pi}\int_{|t| \le \frac{\pi}{2}}
\exp(-itk_n) \xi^n(t) dt
$$
and by substituting $s=t\sqrt n$, this is equal to
$$
= \frac{1}{2\pi}\int_{|s| \le \sqrt n\frac{\pi}{2}}
 \exp(is\frac{k_n}{\sqrt n}) \xi^n(\frac{s}{\sqrt n}) ds
$$
We emphasize that it is fundamental to precisely know the support 
of the values of the variables $X_i$ (i.\ e.\ the minimal lattice
containing these values) since the form of
the inversion formula depends on this. Moreover, the fact
that the $X_i$s take their values on a lattice of span $2$
was used in our first heuristic formula, too.  (This information is
encapsulated in the minimality condition of Theorem 7.)
To prove our desired statement we write
$$
\begin{gathered}
\left| 2\pi\left[ \sqrt nP(W_n=k_n)- 2 \phi(\frac{k_n}{\sqrt n})\right]\right|
 \\ \le
\int_{|s|\le A}\left| \xi^n(\frac{s}{\sqrt n}) - \gamma (s)\right| ds
+ \int_{|s| \ge A} \gamma(s) ds
+ \int_{A \le |s| \le \varepsilon \sqrt n}\left|\xi^n(\frac{s}{\sqrt n})\right|ds
+ \int_{\varepsilon \sqrt n \le |s| \le \sqrt n\frac{\pi}{2}} \left|
\xi^n(\frac{s}
{\sqrt n}) \right|ds\\ = I + II + III + IIII
\end{gathered}
$$
For making the right hand side sufficiently small, we will first select 
$A$ to be sufficiently large and then $\varepsilon$ sufficiently small. Thus
$II$ can be made arbitrarily small, and, for fixed $A$, $I$ will also be small
by the CLT (as a matter of fact, the smallness of $I$ would also
follow from our forthcoming argument for handling $III$. 
In fact, $III$ is a quite
interesting term. By expanding $\xi(t)$ in a power series in the neighbourhood
of $0$, one can easily see that, if $\varepsilon$ is sufficiently small,
then for $|s| \le \varepsilon$ one has
$$
| \xi^n(\frac{s}{\sqrt n}) | \le 1 - \frac{s^2}{4} \le \exp(-\frac{s^2}{4})
$$
As a consequence one obtains
$$
III \le \int_{A \le |s| \le \varepsilon \sqrt n}\exp(-\frac{s^2}{4}) ds
\le \int_{A \le |s|}\exp(-\frac{s^2}{4}) ds
$$
and this term is also small if $A$ is large. Finally, the smallness of $IIII$
follows from the fact that, in the interval
$\varepsilon \sqrt n \le |s| \le \sqrt n\frac{\pi}{2}$,
the term $\left|
\xi(\frac{s}
{\sqrt n}) \right|$ is uniformly bounded away from $1$ from above, and
consequently for the integrand in $IIII$ we have an exponentially 
collapsing upper bound.

\end{document}